\documentclass[12pt]{amsart}
\usepackage{latexsym,amssymb,amsmath,amsthm}
\usepackage{listings}
\usepackage{float}
\usepackage{graphicx}
\usepackage{multicol}
\newtheorem{thm}{Theorem}

\theoremstyle{definition}

\newcommand{\C}{\mathbb{C}}   
\def\O{\mathcal{O}}  \def\G{\mathcal{G}}  \def\P{\mathcal{P}} 
\def\G{\mathcal{G}} \def\B{\mathcal{B}}  \def\C{\mathcal{C}} \def\E{\mathcal{E}}
\def\F{\mathcal{F}} 

 \def\Det{\rm{Det}}  \def\dim{\rm{dim}} 
 \def\deg{{\rm deg}} \def\sign{{\rm sign}}  
 \def\deg{\rm{deg}} \def\str{\rm{str}}
\def\Binomial#1#2{{#1\choose #2}}

\title[Structures within Topological Graph theory]{Classical mathematical structures within topological graph theory}
\author{Oliver Knill}
\date{February 8, 2014}
\address{
        Department of Mathematics \\
        Harvard University \\
        Cambridge, MA, 02138
        }
\subjclass{68-xx, 51-xx }
\keywords{Graph theory, topology, geometry}

\begin{document}
\maketitle
\begin{abstract}
Finite simple graphs are a playground for classical areas of 
mathematics. We illustrate this by looking at some results.
\end{abstract} 

\section{Introduction}

These are slightly enhanced preparation notes for a talk given at the joint AMS meeting of
January 16, 2014 in Baltimore. It is a pleasure to thank the organizers, 
{\bf Jonathan Gross} and {\bf Tom Tucker} for the invitation to participate
at the special section in topological graph theory. \\

A first goal of these notes is to collect some results which hold unconditionally
for any {\bf finite simple graph} without adding more structure.
Interesting are also results which hold for specific classes of graphs like 
{\bf geometric graphs}, graphs of specific dimension $d$ for which the unit spheres 
are $(d-1)$-dimensional geometric graphs. Such graphs behave in many respects like
manifolds. We use that in full generality, any finite simple graph has a natural higher 
simplex structure formed by the presence of complete subgraphs present in the graph. 
Finite simple graphs are intuitive and harbor translations of theorems in mathematics which 
in the continuum need machinery from tensor analysis, functional analysis, complex analysis
or differential topology. Some concepts in mathematics can be exposed free of technicalities 
and could be taught very early on. Here is an example of a result which mirrors a Lefschetz fixed 
point theorem in the continuum and which explains why seeing graphs as higher dimensional objects is 
useful. The result which follows from \cite{brouwergraph} is a discretization  of a fixed point 
result in the continuum:  \\

Given an automorphism $T$ of a triangularization of an $n$-dimensional sphere. Assume $T$ is 
orientation reversing if $n$ is odd and orientation preserving if $n$ is even.
Then $T$ has at least two fixed simplices. The reason is that the Lefschetz number of $T$ is $2$ and
that the sum of the indices of fixed point simplices adds up to $2$. 
For an orientation preserving automorphism of a triangularization of a 2-sphere for example, fixed points are either 
triangles, vertices or edges. An other consequence of \cite{brouwergraph} is that any automorphism of a tree has 
either a fixed vertex or edge, a result which follows also from \cite{NowakowskiRival}.

\section{Results}

The following theorems hold unconditionally for any finite simple graph $G=(V,E)$. They are 
well known in the continuum, but are relatively new within graph theory. 
Here $T$ is an automorphism of the graph with Lefschetz number $L(T) = \sum_{k=0} (-1)^k {\rm tr}(T_k)$, where
$T_k$ is the linear map induced on the $k$'th cohomology group $H^k(G)$. Note that $H^k$ is not only 
considered for $k=0,1$ as in some graph theory literature. 
A triangularization of a $n$-dimensional sphere for example
has ${\rm dim}(H^0(G)) = {\rm dim}(H^n(G)) = 1$ and all other $H^k(G)=\{0\}$ leading to the Betti vector 
$\vec{b}=(1,0,\dots,0,1)$. Only for triangle-free graphs, where $v_0=|V|,v_1=|E|$ determine the Euler characteristic
$\chi(G)=v_0-v_1$ we have by Euler-Poincar\'e $\chi(G)=b_0-b_1$, the difference between the number of components and
genus. For us, triangle-free graphs are special geometric objects with dimension $\leq 1$ and dimension $1$ if there
are no isolated vertices. In other words, only for triangle free graphs, we consider graphs as ``curves".
The function $K(x)= \sum_{k=0}^{\infty} (-1)^k V_{k-1}(x)/(k+1)$ is the Euler curvature, where
$V_k(x)$ is the number of $K_{k+1}$ subgraphs of $S(x)$ and $V_{-1}(x)=1$. The integer $i_T(x) = \sign(T|x) (-1)^{k(x)}$ 
is the {\bf degree} of $T$ for the simplex $x$ and $i_f(x)=1-\chi(S(x) \cap \{ y\; |\; f(y)<f(x) \; \}$ 
is the {\bf index} of a function $f$ at a vertex $x$. If $d$ is the exterior derivative matrix, the matrix $L=d d^* + d^* d$ 
is the form-Laplacian which when restricted to $k$-forms is denoted $L_k$. For a complete graph $G=K_{n+1}$ for example,
$L$ is a $2^{n} \times 2^{n}$ matrix which decomposes into blocks $L_k$ of $B(n,k) \times B(n,k)$ matrices, where 
$B(n,k)=n!/(k!(n-k)!)$. An injective function $f$ on $V$ has a critical point $x$, if $S^-(x) = \{ y  \in S(x) \; | \; f(y)<f(x) \}$
is not contractible. Let ${\rm crit}(G)$ the maximal number of critical points, an injective function $f$ on the vertices can have. 
Let ${\rm tcap}$ denote the minimal number of in $G$ contractible graphs which cover $G$. A graph is {\bf contractible} if
a sequence of homotopy steps, consisting either of pyramid extensions or removals, brings it down to the one point graph
$K_1$. Let ${\rm cup}(G)$ be the {\bf cup length} of $G$ as used in the discrete in \cite{josellisknill}. We compute the 
cup length in an example given in the section with remarks. 
The minimal number $m$ of $k$-forms $f_j$ with $k \geq 1$ in this algebra with the property that
$f_1 \wedge f_2 \cdots \wedge f_k$ is not zero in $H^m(G)$ is called the {\bf cup length}. It is 
an algebraic invariant of the graph. 

\begin{thm}[Gauss-Bonnet] $\sum_{x \in V} K(x) = \chi(G)$. \label{1} \end{thm}
\begin{thm}[Poincar\'e-Hopf] $\sum_{x \in V} i_f(x) = \chi(G)$. \label{2} \end{thm}
\begin{thm}[Index expectation] ${\rm E}[ i_f(x)]= K(x)$. \label{3} \end{thm}
\begin{thm}[Lefschetz] $L(T) = \sum_{T(x)=x} i_T(x)$. \label{4} \end{thm}
\begin{thm}[Brower]$G$ contractible, then $T$ has a fixed simplex. \label{5}\end{thm}
\begin{thm}[McKean-Singer] $\str(e^{-t L}) = \chi(G)$. \label{6} \end{thm}
\begin{thm}[Hodge-DeRham] ${\rm dim}({\rm ker}(L_k)) = {\dim}{H^k(G)} = b_k$ \label{7} \end{thm}
\begin{thm}[Ljusternik-Schnirelmann] ${\rm cup}(G) \leq {\rm tcap}(G) \leq {\rm crit}(G)$. \label{8} \end{thm}

For Theorem~(\ref{1}) see \cite{elemente11,cherngaussbonnet}, for Theorem~(\ref{2}) see \cite{poincarehopf},
for Theorem~(\ref{3}) see \cite{indexexpectation}, for Theorem~(\ref{4}) see \cite{brouwergraph}. In \cite{KnillTopology}
is a Kakutani version. For Theorem~(\ref{5}), see \cite{knillmckeansinger} where also Theorem~(\ref{7}) appears.
As shown in the appendix of that paper, the result is very close to the classical result
\cite{McKeanSinger}.
Theorem~(\ref{8}) in \cite{josellisknill} is especially striking because it relates an {\bf algebraic quantity} 
${\rm cup}$ with a {\bf topological quantity} ${\rm tcap}$ and an {\bf analytic quantity} ${\rm crit}$. The discrete result
\cite{josellisknill} is identical to the continuum result, in the continuum however often other counting 
conditions are used. We can make both ${\rm tcap}$ and ${\rm crit}$
homotopy invariant by minimizing over all graphs homotopic to $G$. Then the two inequalities
relate three homotopy invariants: and algebraic, a topological and an analytical one. \\

The following theorems were known for graphs already, sometimes in other incarnations. We assume the graph to be connected.
For the statement of Kirchoff's theorem we avoid pseudo determinants by using a {\bf Google damping matrix} 
$P_{ij} = 1/n$ which when added just shifts the eigenvalue $0$ to $1$ so that $\Det(L)=\det(P+L)$, where
$\Det$ is the pseudo determinant \cite{cauchybinet}, 
the product of nonzero eigenvalues. Lets write shortly $r$-forest for {\bf rooted spanning forests}.
For the Riemann-Hurwitz statement, we assume that $G/A$ is a graph and $A$ a subgroup of order $n$ of the 
Automorphism group of the graph. 

\begin{thm}[Ivashchenko] Cohomology is a homotopy invariant. \end{thm}
\begin{thm}[Kirchoff] $\det(P+L)/n$ is the $\#$ of maximal trees. \end{thm} 
\begin{thm}[Chebotarev-Shamis] $\det(1+L)$ is the $\#$ of r-forests. \end{thm}
\begin{thm}[Stokes-Gauss] $\sum_{x \in \G} df(x) = \sum_{x \in \delta \G} f(x)$. \end{thm}
\begin{thm}[Euler-Poincar\'e] $\sum_{k=0}^{\infty} (-1)^k \dim(H^k(G)) = \sum_{k=0}^{\infty} (-1)^k v_k$. \end{thm}
\begin{thm}[Riemann-Roch] $r(D) - r(K-D) = \chi(G) + \deg(D)$. \end{thm}
\begin{thm}[Riemann-Hurwitz] $\chi(G) = n \chi(G/A) - \sum_{x \in \G} (e_x-1)$. \end{thm}

For a short exposition on Stokes, see \cite{knillcalculus}. 
Ivashcheko's result \cite{I94} is of great importance because it allows us to take a large complex network, homotopy 
deform it to something smaller, then compute the cohomology using linear algebra.
Better even than homotopy shrinking procedures is a \v{C}ech approach: find a suitable topology on the graph \cite{KnillTopology}
and compute the cohomology of the nerve graph. Kirkchoff's theorem gives the order of the Jacobian group of 
a graph and was primarily the reason we got more interested in it \cite{cauchybinet}. 
For Chebotarev-Shamis \cite{ChebotarevShamis1}, there is an elegant proof in \cite{knillforest} using
classical multilinear algebra. 
For Riemann-Roch, see \cite{BakerNorine2007}, for Riemann-Hurwitz \cite{TuckerKnill} we can assure
that $G/A$ is a graph if we take {\bf simple group actions} which prevent that the quotient graph 
gains higher dimensional simplices. Riemann-Hurwitz holds for pretty arbitrary graphs. The ramification
points can be higher-dimensional simplices so that the result holds for pretty general group actions on 
graphs. Stokes holds more generally for {\bf chains} as Poincar\'e knew already. It is
part of graph theory if $\delta G$ is a graph. Assume $f$ is a $k$-form. Summing over 
the set of simplices means that we sum over all $k$-dimensional simplices in $\G$ or its boundary
$\delta \G$, which is the set of simplices in $\delta G$. 
Stokes can be abbreviated as $\langle G,df \rangle = \langle \delta G, f \rangle$, indicating
that the {\bf exterior derivative} $d$ is dual to the {\bf boundary operation} $\delta$. 
\cite{BakerNorine} have formulated a Riemann-Roch theorem for $1$-dimensional multi-graphs.
It is formulated here for triangle-free graphs to indicate that we neglect the higher-dimensional
structure. The {\bf principal divisor} $K$ is $-2$ times the curvature indicated that Riemann-Roch
is related to Gauss-Bonnet, but the result is definitely deeper.
Divisors are integer-valued functions on vertices. As in the continuum, Riemann-Roch turns out
to be a sophisticated Euler-Poincar\'e formula relating analytically and combinatorically defined 
quantities. A higher-dimensional version will need a discrete analogue of sheaf cohomology.  \\

Riemann-Hurwitz holds in full generality for chains
and just reflects the Burnside lemma for each dimension \cite{TuckerKnill}. For the theorem to work within graph theory
we have to insist that the "orbifold" $G/A$ is a graph. This is not always the case as already a $Z_n$ action on $C_n$
shows, for which the quotient is no more a finite simple graph. Riemann-Hurwitz has appeared in graph 
theory before, but only in the case when graphs are looked at as discrete analogues of 
algebraic curves or Riemann surfaces see \cite{MednykhMednykh} (a reference, I owe Roman Nedela). 
Assume a finite group $A$ of order $n$ acts on a finite graph $G$ by automorphisms.
The theorem sees $G$ as a ramified cover of the chain $H=G/A$ and ramification indices $e_x-1$
with $e_x=1+\sum_{a \neq 1, a(x)=x} (-1)^{k(x)}$, where $k(x)$ is the dimension of the simplex $x$. \\

Lets call a function $f$ on the vertices a {\bf Morse function} if adding a new point along the 
filtration defined by $f$ changes none or exactly one entry $b_m$ in the {\bf Betti vector} 
$\vec{b}=(b_0,b_1, \dots )$.  If the entry $b_m$ is increased, this corresponds to add a $m$-dimensional 
``handle" when adding the vertex.  Now, the index $i_f(x)$ of 
each critical point is by definition either $1$ or $-1$. Adding a zero-dimensional handle for example increases $b_0$ and
also the number of connected components, adding a one-dimensional handle increases $b_1$ and has the effect of 
``closing a loop". One can write $i_f(x)=(-1)^{m(x)}$ where $m(x)$ is the {\bf Morse index},
which is the integer $m \geq 0$ at which the Betti number $b_m$ has changed. In the continuum, $m(x)$ is the dimension of the
stable manifold at the critical point. For the minimum, $m(x)=0$ and for a maximum $m(x)=n$, where $n$ is the dimension. 
Denote by $c_m$ the number of critical points of Morse index $m$. Then

\begin{thm}[Weak Morse] $\chi(G) = \sum_k (-1)^k c_k$ and $b_m \leq c_m$. \end{thm}
\begin{thm}[Strong Morse] $\sum_{k=0}^{m} (-1)^k b_{m-k} \leq \sum_{k=0}^m (-1)^k c_{m-k}$. \end{thm}

The proof is by induction by adding more and more points to the graph. The definition of Morse function
has been made in such a way that the inequalities remain true under the induction step of adding an other vertex.
The induction foundation holds because the results hold for a one point graph $K_1$. Discrete Morse theory has been 
pioneered in a different way \cite{Forman1999,Forman2003} and is more developed. \\

Discrete PDE dynamics reduces to almost trivial linear algebra in the graph theoretical setup if the PDE 
is linear and involves the Laplacian on the geometry as many problems do. We mention this because in the 
continuum, there is a relatively large technical overhead with integral operators, 
distributions and functional analysis just to make sense of objects like Greens functions. One reason is
that in the continuum, the involved operators are unbounded, making the use of functional analysis unavoidable
for example just to be able to establish spectral properties like elliptic regularity or even establish the
existence of solutions of the PDEs in suitable function spaces. 
The matrix $L_h$ is the form-Laplacian restricted to the complement of the kernel. It operates on general
$k$-forms, and these matrices are easy to write down for any finite simple graph. 
We could write $d_0={\rm grad}, d_0^* = {\rm div}, d_0^* d_0 = L_0 = \Delta, d_1 = {\rm curl}$.
Electromagnetic waves, heat, gravity all make sense on a general finite simple graph. 
The {\bf Newton gravitational potential} $V$ for example satisfies $L V = \rho$, where $\rho$ is the 
{\bf mass density}, a function on the vertices. We can get from the mass density the {\bf gravitational field}
$F=dV$, a one-form. Note that this works on any finite simple graph: the Laplacian defines a natural 
Newton potential. For electromagnetism, we get from the {\bf charge-current one-form} $j$
the {\bf electromagnetic potential} $A$ and from that the {\bf electromagnetic field} $F=dA$. 
When studiing the {\bf wave equation} it becomes apparent how useful it is to see graphs as discrete 
Riemannian manifolds. While the Dirac operator $D=d+d^*$ is a cumbersome object in the continuum, it is natural in 
the discrete as it leads immediately to a basis for the cohomology groups $H^k(G)$ - which are vector spaces -
by computing the kernel of each block matrix $L_k$ of $L=D^2 = d d^* + d^* d$. The matrix $D$ is the crux of
the story because it encodes the full exterior derivative $d$ in the upper triangular part, its adjoint $d^*$ in 
the lower triangular part.  Evolving partial differential equations on a graph reduces to matrix exponentiation in 
linear algebra. Gravity lives on zero-forms, electromagnetism on one-forms, the weak force on two-forms and the 
strong force on three-forms. We mention this because in the continuum, there are chapters of books dedicated 
to the problem just to find the electromagnetic field to a current and charge distribution. To find the 
gravitational field $F$ of a mass distribution $\rho$ on a Riemannian manifold, we have to compute 
the {\bf Green's kernel} which already uses the language of distributions in the continuum. The Poisson 
equation is a system of linear equations. The heat equation is $u'=-Lu$, the wave equation is $u''=-Lu$
the Maxwell equations are $dF=0,d^*F=j$ for a one form $j$. A Coulomb gauge $d^*A=0$ reduces Maxwell to 
a Poisson equation for $1$-forms, as in the continuum. 

\begin{thm}[Fourier] $e^{ -Lt} u(0) $ solves the heat equation $u' = -L u$. \end{thm}
\begin{thm}[d'Alembert] $\cos(Dt) u_0+ \sin(Dt) D^{-1}u'_0$ solves wave. \end{thm}
\begin{thm}[Poisson] $L_h^{-1} g$ solves the Poisson equation $Lu = g$  \end{thm} 
\begin{thm}[Maxwell] $A = L_h^{-1} j, F=dA$ solves Maxwell. \end{thm} 
\begin{thm}[Gauss] $d^* F = \rho$ defines gravity $F=dV$ by $V=L_h^{-1} \rho$. \end{thm}
\begin{thm}[Hopf-Rynov] For $x,y \in V$, exists $v$ with $\exp_x(v) =y$. \end{thm}
\begin{thm}[Toda-Lax] $D'=[B,D]$ with $B(0)=d-d^*$ is integrable. \end{thm}

The heat equation is important because we could use it to find harmonic $k$-forms if it were not just
given already as an eigenvalue problem. 
The wave equation $(d^2/dt^2 +L) = 0$ can be factored $(d/dt + i D) (d/dt -iD) \psi=0$ leading
to {\bf Schr\"odinger equations} $\psi'= \pm i D \psi$ for the Dirac operator $D$,and complex quantum 
wave $\psi(t)=u(t) + i D^{-1} u'(t)$ encoding position $u(t) = {\rm Re}(\psi(t))$ and velocity
$u'(t) = D {\rm Im}(\psi(t))$ of the classical wave. We write $\exp_x(v) = \psi(t)$ if $\psi(0) = x+D^{-1} v$. 
It is convenient for example to use the wave flow as a discrete analogue of the {\bf geodesic flow} which in 
physical contexts has always been given by light evolution: we measure distances with light.
What happens in the discrete is that Hopf-Rynov only can be realized
by looking at a quantum dynamical frame work. It is absolutely futile to try to find a notion of 
tangent space and exponential map on a discrete level by only using paths on the graph. The reason is that
for any pair of points $x,y$ we want to have an element $v$ in the tangent space so that $exp_x(v)=y$. If 
the graph is large then the tangent space needs to be large. One can try to look at equivalence classes of paths
through a vertex but things are just not working naturally. Linear algebra is not only easier, it is also 
how nature has implemented the geodesic flow: as motions of particles satisfying quantum mechanical rules. 
Hopf-Rynov in the graph case is now very easy: for any two vertices $x,y$ there is
a unique initial velocity of a wave localized initially at $x$ so that at a later time $T$ it is at $y$. This is just
linear algebra. As in the continuum there is an uncertainty principle: while we can establish to have a particle passing through
two points (``knowing the velocity") prevents us to know the position exactly at other times. But we can have a classical
motion on the graph as the vertex on which the probability density of the particle position is highest. In the graph
case, the quantum unitary evolution happens in a compact unitary group so that there is always Poincar\'e recurrence
evenso as in classical mechanics, the return times are huge already for relatively small graphs. 
We can use the wave equation naturally to measure distances, which does this much better than the
very naive geodesic distance. Naturally, different types of particles - as waves of differential forms - travel with different
velocities too. The point of view that a graph alone without further input allows to study relatively complex physics
has been mentioned in \cite{eveneuler,IsospectralDirac}. We don't even need to establish initial conditions since the isospectral 
deformation of $D(t)$ does that already. Which isospectral deformation do we chose? It turns out not to matter. The physics
is similar. A finite simple graph leads to physics without further assumptions. The only input is the graph. 
The details of the dynamical system looks difficult: 
find relations between the speeds with which different discrete differential forms $f \in \Omega_k$ move. 
Since the Dirac operator links different forms, this can not be studied on each $k$-form sector 
$\Omega_k$ separately. Only the Laplacians $L_k$ leave those sectors invariant - sectors which
informally could be thought of as gravitational, electric, weak and strong. The symmetry translations
given by isospectral deformations $D(t)$ allows to 
study the column vectors of $D(t)$ which asymptotically solve the wave equation
and especially explore the geometry on $\Omega= \oplus \Omega_k$ by computing
distances between various simplices. This leads to speed relations between various particles. Even for smaller
graphs, the dynamics is complex without the need for any further input. 
The dynamics can be simulated on the computer. Evenso the picture is naive, it looks like 
a wonderful playground for experimentation. In trying to figure out, which graphs are natural we have looked in
\cite{eveneuler} at the Euler characteristic is a natural functional.
The isospectral deformation of the Dirac operator $D=d+d^*$ can be done in such a way that it 
becomes in the limit a wave evolution. It is interesting that this modification of the flow, where we take
$D(t) = d+d^*+b$. The deformation of the Dirac operator studied in \cite{IsospectralDirac} 
can be modified so that it becomes complex. We have then more symmetry with a selfadjoint
operator $D=d+d^*+b$ and an antisymmetric operator $B=d-d^*+ib$. 
The dynamics (and geometry given by the differential forms $d$) 
now becomes complex even if we start with a real structure. With or without the modification, it 
provides an isospectral deformation of $D$. The Laplacian $L$ stays the same 
so that the geometric evolution can not be seen on a classical level, except when looking at the d'Alembert solution
of the wave equation, where $D$ enteres in the initial condition. The symmetry group of isospectral Dirac operators on 
a  graph provides a natural mechanism for explaining why geometry expands with an inflationary start. See also \cite{eveneuler}. 
This is not due to some special choice of the deformation but is
true for any deformation starting with a Dirac operator which has no diagonal part. The deformed operator 
will have some diagonal part which leads to "dark matter" which is not as geometric as the side diagonal part. 
This simple geometric evolution system \cite{IsospectralDirac2,IsospectralDirac}
allows to deform both Riemannian manifolds for which the deformed $d$ are pseudo differential operators.
It is an exciting system and can not be more natural because one is forced to consider it when taking
quantum symmetries of a graph seriously. Besides establishing basic properties like universal expansion with inflationary 
start and camouflaged supersymmetry, establishing the limiting properties of the system has not yet been done. \\

A graph has {\bf positive curvature} if all 
sectional curvatures are positive. In the following, a {\bf positive curvature graph} is a 
{\bf geometric positive curvature graph}, meaning that there is $d$ such that 
for all vertices $x$ the sphere $S(x)$ is a $d-1$ dimensional geometric graph which is a sphere
in the sense that the  minimal number of critical points, an injective function on the graph 
$S(x)$ can take is $2$. The set of all $d$-dimensional positive curvature graphs is called $\P_d$. 
A {\bf sectional curvature} is the curvature of an embedded wheel graph.
We denote by $K(x)$ the Euler curvature, with ${\rm Aut}(G)$ the group of automorphisms and with 
${\rm Aut}_+(G)$ the group of orientation-preserving automorphisms of a geometric graph $G$ and
with $\F$ the set of fixed simplices of $T$. 

\begin{thm}[Flatness] Odd dimensional graphs have $K(x) \equiv 0$. \label{10} \end{thm}
\begin{thm}[Bonnet] Positive curvature graphs have diameter $\leq 3$.  \end{thm}
\begin{thm}[Synge] $2m$-dim positive curvature $\Rightarrow$ simply connected.\end{thm}
\begin{thm}[Bishop-Goldberg] $G \in \P_{2d}$ $\Rightarrow$ $b_0=b_3=1,b_1=b_2=0$. \end{thm}
\begin{thm}[Weinstein] In $G \in \P_4$, $T \in {\rm Aut_+}(G)$ has $|\F| \geq 2$. \end{thm}

For Theorem~(\ref{10}), see \cite{indexformula}.  \\

A connection with different fields of mathematics comes 
through zeta functions $\zeta(s) = \sum_{\lambda>0} \lambda^{-s}$, where $\lambda$ runs
over all positive eigenvalues of the Dirac operator $D$ of a graph $G$.
These entire functions leads to connections to basic complex analysis. Why using the Dirac operator? One
reason is that for circular graphs $G=C_n$, the zeta function has relations with the corresponding
zeta function of the Dirac operator $i \partial_x$ of the circle $M=T^1$, for which the zeta function is
the {\bf standard Riemann zeta function}. Also, we have $\zeta(2s)$ as the zeta function of the Laplace operator.
It is better to start with $\zeta(s)$ and go to $\zeta(2s)$ than looking at the Laplace zeta function $\tilde{\zeta}(s)$
and then have to chose branches when the square root $\tilde{\zeta}(s/2)$.
Again, many technicalities are gone since we deal with entire functions so that there is no surprise that for 
discrete circles $C_n$. While there is no relation with the Riemann hypothesis 
- the analogue question for the circle $S^1$ - it is still interesting because the proof \cite{KnillZeta}
has relations with single variable calculus.

\begin{thm}[Baby Riemann] Roots of $\zeta_{C_n}(s)$ $\rightarrow$ ${\rm Re}(s) = 1/2$. \end{thm}

Finally, lets mention an intriguing class of {\bf orbital networks} which we discovered first together
with Montasser Ghachem and which has much affinities to realistic classical 
networks and leads to questions of number theoretical nature. We consider $k$ polynomials maps $T_k$ on the ring
$V=Z_m$ or multiplicative group $Z_m^*$ 
which we take as vertices of a graph. Two vertices $x,y$ are connected if there is $k$ such that $T_k(x)=y$ or
$T_k(y)=x$. Here are ``miniature example results" \cite{KnillOrbital1,KnillOrbital2,KnillOrbital3}:

\begin{thm} $(Z_n^*,2x)$ is connected iff $n=2^m$ or $2$ is primitive root. \end{thm}
\begin{thm} $(Z_n^*,x^2)$ is connected iff $n=2$ or $n$ is Fermat prime. \end{thm}
\begin{thm} $(Z_n,\{ 2x,3x+1 \})$ has $4$ triangles for prime $n>17$. \end{thm}
\begin{thm} $(Z_n^*,\{ x^2,x^3 \})$ is connected iff $n$ is Pierpont prime. \end{thm}
\begin{thm} $(Z_n,T)$ has no $K_4$ graphs and $\chi \geq 0$. \end{thm}
\begin{thm} $(Z_p,x^2+a,x^2+b)$ with $a \neq b$ has $\chi<0$ for large primes. \end{thm}
\begin{thm} $(Z_p,x^2+a)$ has zero, one or two triangles. \end{thm}

More questions are open than answered: we still did not find a Colatz graph $(Z_n,2x,3x+1)$ 
which is not connected. We have found only one example of a quadratic orital network with 3 different generators
which is not connected: the only case found so far $(Z_{311},x^2+57,x^2+58,x^2+213)$ and checked until $p=599$
but since billions of graphs have to be tested for connectivity, this computation has slowed down considerably, 
now using days just for dealing with one prime. We als see that for primes $p>23$ and $p \leq 1223$ all 
quadratic orbital networks with two different generators are not planar. Also here, we quickly reach computational
difficulties for larger $p$. These are serious in the sense that the computer algebra system refuses to decide
whether the graph is planar or not. 

\section{Remarks}

Having tried to keep the previous section short and the statements concise, the following remarks are
less polished. The notes were used for preparation similarly as \cite{DiracKnill} for the more linear algebra
related issues. \\

{\bf Dimension.}
Classically, the {\bf Hausdorff-Uhryson inductive dimension} of a topological space $X$ is defined
as ${\rm ind}(\emptyset)=-1$ and ${\rm ind}(X)$ as the smallest $n$ such that for every 
$x \in X$ and every open set $U$
of $x$, there exists an open $V \subset U$ such that the boundary of $V$ has
dimension $n-1$. For a finite metric space and a graph in particular, 
the inductive dimension is $0$ because every singleton set $\{x \;\}$ is open.
It can become interesting for graphs when modified: 
dimension \cite{elemente11} for graphs is defined as
${\rm dim}(\emptyset)= -1$, ${\rm dim}(G) = 1+\frac{1}{|V|} \sum_{v \in V} {\rm dim}(S(v))$,
where $S(v)$ is the unit sphere. Dimension associates a rational number to every vertex $x$, 
which is equal to the dimension of the sphere graph $S(x)$ at $x$. 
The dimension of the graph itself is the average over the dimensions for all vertices $x$.
We can compute the expected value $d_{n}(p)$ of the dimension on Erdoes-R\'enyi probability spaces 
$G(n,p)$ recursively: $d_{n+1}(p) = 1+\sum_{k=0}^n \Binomial{n}{k} p^k (1-p)^{n-k} d_k(p$, 
where $d_0=-1$. Each $d_n$ is a polynomial in $p$ of degree $\Binomial{n}{2}$. \cite{randomgraph}
Triangularizations of $d$-dimensional manifolds for example have dimension $d$. \\

{\bf Homotopy.}
Classically, two topological spaces $X,Y$ are homotopic, if there exist continuous maps $f:X \to Y$,
$g:Y \to X$ and $F:X \times [0,1] \to X, G:Y \times [0,1] \to Y$ such that $F(x,0)=x,F(x,1)=g(f(x))$
and $G(y,0)=y, G(y,1) = f(g(y))$. For graphs, there is a notion which looks different at first: it 
has been defined by Ivashchenko \cite{I94} and was refined in \cite{CYY}. 
We came to this in work with Frank Josellis via Morse theory as we worked on \cite{poincarehopf}: it became
clear that the set $S^-(x) = \{ y \; | \; f(y)<f(x) \; \}$ better has to be contractible if nothing 
interesting geometrically happens. Indeed, we called in \cite{poincarehopf} the quantity 
$1-\chi(S^-(x))$ the index. Given a function $f$ 
on a graph, we can build a Morse filtration $\{ f \leq c \}$ and keep track of the moment, when something
interesting happens to the topology. Values for which such a thing happens are {\bf critical values}. 
If we keep track of the Euler characteristic, using the fact that unit balls have Euler characteristic $1$
and $\chi(A \cap B) = \chi(A) + \chi(B) - \chi(A \cap B)$ by counting, we immediately get $\sum_x i_f(x) = \chi(G)$
which is the Poincar\'e-Hopf theorem. We can look at extension steps $X \to Y$, in which a new vertex $z$
is attached to a {\bf contractible} part of the graph. This is a homotopy step and is defined recursively because
{\bf contractible} means homotopic to $K_1$. The reverse step is to take a 
vertex $z$ with contractible unit sphere and remove both $z$ and all connections to $z$. A graph $G$ 
for which we can apply a finite sequence of such steps is called a {\bf homotopy}. Homotopy preserves
cohomology and Euler characteristic. But as in the classical case, 
homotopy does not preserve dimension: all complete graphs $K_n$ are homotopic. \\

{\bf Continuity.}
While natural notions of homotopy and cohomology exist for graphs, we need a notion of 
homeomorphism which allows to deform graphs in a rubber geometry way. We want dimension, homotopy and cohomology
and connectivity as in the continuum. Any notion of topology on a graph which fails one of these conditions is
deficient. Topology is also needed for discrete sheave constructions.
Classically, a notion of homeomorphism allows to do subdivisions of edges of a graph to get a homeomorphic graph. 
This notion is crucial for the classical Kuratowski theorem which tells that a graph is not planar if and only it 
contains a homeomorphic copy of either $K_5$ or $K_{3,3}$. 
This classical notion however does not preserve dimension, nor cohomology, nor Euler characteristic in general: 
the triangle is homotopic to a point, has Euler characteristic $1$ and trivial cohomology. After a subdivision of the 
vertices we end up with the graph $C_6$ which has Euler characteristic $0$, is not homotopic to a point and has
the betti vector $b_1=1$. In \cite{KnillTopology}, we have introduced a notion of homeomorphism which 
fulfills everything we wish for. The definition: 
a {\bf topological graph} $(G=(V,E),\B,\O)$ is a standard topology $\O$ on $V$ generated by a subbasis $\B$ of 
contractible sets which have the property that if two basis elements  $A,B$ intersect under the 
{\bf dimension condition} $\dim(A \cap B) \geq {\rm min}(\dim(A),\dim(B))$ then 
the intersection is contractible. The {\bf nerve graph} $(\B,\E)$ is given by the set $\E$ consisting of all
pairs $A,B$ for which the dimension assumption holds. This nerve graph is asked to be homotopic to the 
graph. A {\bf homeomorphism} of two-topological graphs is a classical homeomorphism between the topologies
so that the subbasis elements correspond and such that dimension is left constant on $\B$. Two graphs $G,H$
are homeomorphic if there are topologies $(\B,\O)$ and $(\C,\P)$ on each
for which they are homeomorphic. Homeomorphic graphs are homotopic,
have the same cohomology and Euler characteristic. Lets call a topological graph $(G,\B,\O)$ {\bf connected} if
$\B$ can not be split into two disjoint sets $\B_1 \cup \B_2$ such that for every $B_i \in \B_i$ the
sets $B_1 \cap B_2$ are empty. For any topology this notion of connectivity is equivalent to path connectivity.
(Also in classical point set topology, a space $X$ is connected if and only if there is a subbasis of 
connected sets for which the nerve graph is connected.)
Every graph has a topology for which the nerve graph is the graph itself. 
With the just given definition of topology, any two cyclic graphs
$C_n,C_m$ are homeomorphic if $n,m \geq 4$. Also, the octahedron and icosahedron are 
homeomorphic. Contractible graphs can 
carry weak topologies so that all trees are homeomorphic with respect to this indiscrete topology.
Graphs are isomorphic as graphs if we take the discrete topology on it, the topology generated by star graphs
centered at wedges. Under the fine topology, two graphs are homeomorphic if and only if they are isomorphic.
As in the continuum, natural topologies are neither discrete nor indiscrete. For 
triangularizations of manifolds, we can take $\B$ to consist of discs, geometric graphs of the same dimension
than the graph for which the boundary is a Reeb sphere. \\

{\bf Trees and forests.}
The matrix tree theorem of Kirkhoff is important because the number of spanning tree is the dimension of the 
Jacobian group of a graph. Higher dimensional versions of this theorem will be relevant in higher dimensions. 
The {\bf Chebotarev-Shamis forest theorem} follows from a general Cauchy-Binet identity,
if $F,G$ are two $n \times m$ matrices, then
$\det(1+x F^T G) = \sum_P x^{|P|} \det(F_P) \det(G_P)$,
where the sum is over all minors \cite{cauchybinet}.  This implies
$ \det(1+F^T G) = \sum_P \det(F_P) \det(G_P)$ and especially a cool {\bf Pythagorean identity}
$\det(1+A^T A) = \sum_P \det^2(A_P)$ which is true for any $n \times m$ matrix. See
\cite{ChebotarevShamis1}. We are not aware that even that special Pythagorean identity 
has appeared already. Note that the sum is over all square sub patterns $P$ of the matrix. \\

{\bf Curvature.}
For a general finite simple graph, the {\bf curvature} is defined as
$K(x)= \sum_{k=0}^{\infty} (-1)^k V_{k-1}(x)/(k+1)$, where
$V_k(x)$ be the number of $K_{k+1}$ subgraphs of $S(x)$ and $V_{-1}(x)=1$. It leads
to the {\bf Gauss-Bonnet-Chern} theorem $\sum_{x \in V} K(x) = \chi(G)$ \cite{cherngaussbonnet}.
These results are true for any finite simple graph \cite{knillcalculus}. 
Gauss-Bonnet-Chern in the continuum tells that for a compact even dimensional 
Riemannian manifold $M$, there is a local function $K$ called Euler curvature
which when integrated produces the Euler characteristic. 
Other curvatures are useful. In Riemannian geometry, one has sectional curvature, which when known along
all possible planes through a point allows to get the Euler curvature, Ricci curvature or scalar curvature. 
In the discrete, the simplest version of sectional curvature is the curvature of an embedded wheel graph. 
if the center of the later has degree $d$, then the curvature is $1-d/6$. This means that wheel graphs 
with more than $6$ spikes have negative curvature and wheel graphs with less than $6$ spikes have 
positive curvature. It is easy to see that any graph for which all sectional curvatures are positive must 
be a finite graph. More generally, if one assumes that positive curvature graphs have positive density 
on any two dimensional geometric subgraph, then one can give bounds of the diameter of the graph similarly
as in the continuum.  Unlike as in differential geometry, curvature appears to be a pretty rigid notion. 
This is not the case: as indicated below, we have proven that curvature is the average over all indices $i_f(x)$ 
when we average over a probability space of functions. This is an integral geometric point of view which 
allows deformation. If we deform the functions $f$, for example using a discrete partial differential equation, 
then the probability measure changes and we obtain different curvatures which still satisfy Gauss-Bonnet
\cite{IsospectralDirac2,IsospectralDirac}. This holds also in the continuum: any probability measure on 
Morse functions on a Riemannian manifold defines a curvature by taking the index expectation on that
probability space. For manifolds embedded in a larger Euclidean space (no restriction of generality by Nash),
one can take a compact space of linear functions in the Ambient space and integrate over the natural 
measure on the sphere to get Euler curvature. What is still missing is an {\bf intrinsic measure} on the 
space of Morse functions of a Riemannian manifold.  \\

{\bf Morse filtration.}
A function $f$ on the vertex set $V$ of a finite simple graph defines a filtration of the graph into 
subgraphs $\emptyset = G_0 \subset G_1 \dots \subset G_n=G$, where each $G_k = \{ v \; | \; f(v) \leq c_k \}$
is obtained from $G_{k-1}$ by adding a new vertex for which the value of $f$ is the next bigger number in the
range of $f$.  The sphere $S(x_k)$ in $G_{k}$ is $S^-(x_k)$. Define
the index of $x$ to be $i_f(x) = 1-\chi(S^-(x))$. We see that $\chi(G_k) - \chi(G_{k-1}) =  i_f(x_k)$ so 
that $\sum_{x \in V} i_f(x) = \chi(G)$. This is the Poincar\'e-Hopf formula for the gradient field of $f$. 
As a side remark, me mention that also the {\bf Poincar\'e-Hopf theorem} can be considered for more 
general spaces. Assume we have a compact metric space $(X,d)$ for which Euler characteristic is defined and for which 
small spheres $S_r(x)$ are homeomorphic for small enough positive $r$. Given a function $f$, 
we call $x$ a {\bf critical point} if for arbitrary small positive $r$, the set  $S_r^-(x) = S_r(x) \cap \{ f(y) < f(x) \}$ 
is not contractible or empty and the index $i_f(x)$ at a critical point $1-\chi(S_r^-(x))$ exists for small enough $r$
and stays the same. In that case, the Poincar\'e-Hopf formula $\sum_{x} i_f(x) = \chi(X)$ holds. \\ 

{\bf Integral geometry.}
There is a great deal of integral geometry possible on finite simple graphs. 
If we take a probability measure on the space of functions like taking $f$ with 
the uniform distribution $[-1,1]$ independent on each node, then the expectation of $i_f(x)$ is curvature. 
This also works for discrete measures like the uniform distribution on all permutations of $\{1, \dots, n\}$. 
By changing the probability measure, we can get different curvatures, which still satisfy the Gauss-Bonnet
theorem. This principle is general and holds for Riemannian manifolds, even in singular cases like
{\bf polytopes} where the curvature is concentrated on finitely many points. 
Lets look at a triangle $ABC$ in the plane and the probability space of 
linear maps $\{ f(x,y) = x \cos(t) + b \sin(t) \; | \; 0 \leq t <2\pi \}$ on the plane which induces functions
on the triangle. Define the index $i_f(x) = \lim_{r \to 0} \chi( S_r(x) \cap \{y \; | \;  f(y) < f(x) \} )$.
The probability to have a positive index away from the corner is zero. This implies that the curvature has its support
on the vertices. The probability to have index $1$ at $A$ is the curvature $K(A)=1 - \alpha/\pi$. 
The sum of all curvatures is $2$. That the sum of the angles $\alpha+\beta+\gamma$ 
in a triangle is equal to $\pi$ can be seen in an integral geometric way as an index averaging result. The same 
is true for polyhedra \cite{Banchoff67} and so by a limiting procedure for general Riemannian manifolds. We have 
not yet found a natural {\bf intrinsic} probability space of Morse functions on a manifold for which the
index expectation is the curvature. Taking an even dimensional Riemannian manifold $M$ with normalized Riemannian volume $m$
we can take the probability space $(M,m)$. 
For any point $x$ one can look at the {\bf heat signature function} $f(y)= K(x,y)$ of the 
{\bf heat kernel} $K(x,y)$. 
The functions $f$ are parametrized by $x$ and are Morse functions. Intuitivly, $f(y)$ is the temperature at a point $y$
of the manifold if the heat source is localized at $x$. The index expectation of this probability space of Morse
function certainly defines a curvature on the manifold for which Gauss-Bonnet-Chern holds. The question is whether it is the
standard Euler curvature. It certainly is so in symmetric situations like constant curvature manifolds. \\

{\bf Divisors.} Lets look at an example for Riemann-Roch. If $G=C_5$ is a cyclic graph then $\chi(G)=1-g=0$ 
and the {\bf principal divisor} is constant $0$. The Riemann Roch theorem now tells $r(D) - r(-D) = {\rm deg}(D)$. 
Given the divisor $D=(a,b,c,d,e)$ we have ${\rm deg}(D)=a+b+c+d+e$. 
The linear system $|D|$ consists of all divisors $E$ for which $D-E$ is equivalent to an effective divisor.
$r(D)$ is  $-1$ if $|D| = \emptyset$ and $r(D) \geq s$ if and only if $|D -E| \neq  \emptyset$ or all 
effective divisors $E$ of degree $s$. If ${\rm deg}(D)=0$, then clearly 
$r(D)=r(-D)=0$. We can assume ${\rm deg}(D)>0$ without loss of generality. Otherwise just change $D$ to $-D$
which flips the sign on both sides. The set of zero divisors is $4$ dimensional. The set of principal divisors
is four dimensional too. The Jacobian group ${\rm Jac}(G)={\rm Div}_0(G)/{\rm Prin}(G)$ has order $5$ and 
agrees with the number of spanning trees in $C_5$ is $5$. It has representatives like $(1,-1,0,0,0)$. 
Given a vertex $x$ the Abel-Jacobi map is $S(v) = [(v)-(x)]$ \cite{BakerNorine2007}. \\

{\bf Geometric graphs}. The emergence of the definition of {\bf polyhedra} and {\bf polytopes} is fascinating.
The struggle is illustrated brilliantly in \cite{lakatos} or \cite{Richeson}. Many authors 
nowadays define polyhedra embedded in some geometric space like \cite{gruenbaum}. 
We have in \cite{elemente11} tried to give a purely graph theoretical definition: a {\bf polyhedron}
is a graph, which can be completed or truncated to become a two dimensional geometric graph. 
The only reason to allow truncations of vertices is so that we can include graphs like the tetrahedron into the class of 
polyhedra. With this notion, all platonic solids are polyhedra: the octahedron and icosahedron are already
two-dimensional geometric graphs, the dodecahedron and cube can be stellated to become a two dimensional 
geometric graph. The notion in higher dimensions is similar. A polytop is a graph which when completed
becomes a $d$-dimensional geometric graph. Now, in higher dimensions, one can still debate whether one
would like to assume that the unit spheres $S(x)$ are all Reeb sphers, or whether one would like to 
weaken this and allow unit sphere to be a discrete torus for example. In some sense this is like allowing 
exotic differential structures on manifolds. \\

{\bf Positive curvature}. 
A classical theorem of Bonnet tells that a Riemannian manifold of positive sectional curvature is compact
and satisfies an upper diameter bound $\pi/\sqrt{k}$. For two-dimensional graphs we can give a list of all the graphs with strictly 
positive curvature and all have diameter $2$ or $3$. We would like to have Schoenberg-Myers which 
needs {\bf Ricci curvature}. The later notion can be defined in the discrete as a function $R(e)$ on edges $e$ given
by the average of the wheel graph curvature containing the edge. The {\bf scalar curvature}
at a vertex $v$ would then just be average of all $R(e)$ with $v \in e$. 
We have not proven yet but expect that only finitely many graphs exist which have strictly positive Ricci curvature
and fixed dimension. Lets go back to positive curvature which implies that every dimensional subgraph has positive curvature.
Positive sectional curvature is a strong assumption in the discrete which prevents to build a discrete positive
curvature projective plane. Since every two dimensional positive curvature graph is a sphere, results are
stronger. Why is there no projective plane of positive curvature? Because the curvature constraint in the discrete produces
too strong pitching. Identifying opposite points of icosahedron
(the graph with minimal possible positive curvature $1/6$ at every vertex) would become higher dimensional.
The statements in the theorems all follow from the following ``geomag lemma":
{\it a closed $2$-geodesic in a positive curvature graph can be extended within $G$ to a two-dimensional
orientable positive curvature graph. } Proof: extend the graph and possibly reuse the same vertices again to build
a two dimensional graph which is locally embedded. The wheel graphs in this $2$-geodesic 2-dimensional subgraph has
4 or 5 spikes due to the positive curvature assumption.
By projecting onto $G$ we see that this is a finite cover of a two-dimensional
embedded graph of diameter $\leq 3$. The later is one of a finite set of graphs which all are orientable
and have Euler characteristic $2$. By Riemann-Hurwitz, the cover must be $1:1$ since there are no two dimensional
connected graphs with Euler characteristic larger than $2$.
Bonnet immediately follows from the Geomag lemma since a geodesic is contained in a two dimensional geodesic surface which
has diameter $\leq 3$. Synge follows also follows from the geomag lemma because all the two dimensional positive curvature
graphs are simply connected and any homotopically nontrivial closed curve can be extended to a two dimensional surface.
Bishop-Goldberg follows now from the Poincar\'e duality.
Weinstein immediately follows from Bishop-Goldberg and the Lefschetz theorem because the Lefschetz number is $2$.
In the continuum, Synge has odd dimensional projective spaces as counter examples for positive curvature
manifolds which are not simply connected. The discrete Bishop-Goldberg
result implies that Euler characteristic is $2$ and so positive. In $d \geq 8$ dimensions this
is unsolved in the continuum but the sectional curvature assumption in the discrete is of course much stronger.
In the continuum, one has first tried to relate positive Euler curvature with positive sectional curvature which fails
in dimensions $6$ and higher \cite{Weinstein71,Geroch}. 
{\bf Ricci curvature} $R(e)$ at an edge $e$ is as the average over all curvatures of wheel graphs which contain $e$.
By extending embedded geodesics to two-dimensional immersed surfaces within $G$, one might get
Schoenberg-Myers type results linking positive Ricci curvature everywhere with diameter bounds.
It could even lead to {\bf sphere theorems} like Rauch-Berger-Klingenberg: Positive curvature graphs are 
triangularizations of spheres. These toy results in the discrete could become more interesting (and get closer to 
the continuum) if the curvature statements are weakened. A graph has {\bf positive curvature of type $(M,\delta$)} if the total
curvature for any geodesic two-dimensional subgraph of diameter $\geq M$ is $\geq \delta$.
Unlike for Riemannian manifolds, where geodesic two-dimensional surfaces in general do not exist, 
the notion of {\bf geodesic surface} is interesting in the discrete: a subgraph $H$ of $G$ is called {\bf L-geodesic} if
for any two vertices $x,y$ with distance $\leq L$ in $H$, the geodesic distance of $x$ and $y$ within $H$ is the same than
the geodesic distance within $G$. Any curve in a graph is a $1$-geodesic and a $2$-geodesic if we shorten each corner in a triangle
by the third side. In other words, a curve is a $2$-geodesic, if we can not do localized homotopy transformations shortening the curve. 
A wheel subgraph is an example of a two-dimensional $2$-geodesic subgraph. In the continuum, geodesic two dimensional surfaces in
$d \geq 3$ dimensional manifolds $M$ do not exist in general: we can form the two-dimensional 
surfaces $\exp_D(x)$, where $D$ is a two-dimensional disc in the tangent space $T_xM$  but for $y$ in this surface the surface
$\exp_D(y)$ intersects $\exp_D(x)$ only in a one dimensional set in general. The geodesic surfaces do not match up. \\

{\bf Structures}. Graph theory allows to illustrate how different fields of mathematics 
like algebra, calculus, analysis, topology, differential geometry, number theory or algebraic geometry interplay.
There are properties of graphs which are of {\bf metric nature} and this reflects through the automorphism
group of the graph consisting of isometries which are graph isomorphisms. {\bf Spectra} of the Laplacian $L$ or Dirac
operator $D$ are examples which are metric properties. 
The spectrum changes under topological deformations or even homotopies. Then there is the group of homeomorphisms
\cite{KnillTopology} of a graph with respect to some topology. This group is in general much richer and more flexible;  
as homeomorphisms do not have to be isometries any more. They still preserve everything we want to be preserved like homotopy,
homology and dimension. Topological symmetries can be weakened when looking at homotopies which provides a much weaker equivalence
relation between graphs. Homotopy still preserves homology and connectivity but does no more preserve dimension. Anything contractible
is homotopically equivalent to the one point graph $K_1$. An even weaker notion is formed by the equivalence classes of graphs with 
the same homology. This algebraic equivalence relation and still weaker than homotopy as we know from the continuum:
there are nonhomotopic graphs which have the same cohomology groups: examples are discrete Poincar\'e spheres. 
There is still room to experiment with {\bf differentiability structures} on a graph.
One possibility for an analogue of a $C^r$-diffeomorphism between two graphs is a homeomorphism which 
has the property that the nerve graph has discs of radius $r$ for which the boundary is a geometric sphere.
One can consider \cite{elemente11} as an attempt to look at curvature in a smoother setup and as the Hopf Umlaufsatz proven there
shows, things are already subtle in {\bf discrete planimetry}, where one looks at regions in the discrete plane. 
Second order curvatures obtained by measuring the size of spheres of radius $1$ and $2$ often fail to satisfy 
Gauss-Bonnet, already in two dimensions. Since it often also works, it is an interesting question to find ``smoothness
conditions" under which the discrete curvature $K(p) = 2 |S_1(p)| - |S_2(p)|$ \cite{elemente11} works for geometric graphs.
The just mentioned higher order curvature could be used as sectional curvature when the spheres are restricted to 
two dimensional subgraphs. The curvature is motivated  from 
Puiseux type formulas $K = \lim_{r \to 0}  \frac{2 |S_{r}| - |S_{2r}|}{2 \pi r^3}$ for 
two dimensional Riemannian surfaces.  \\

{\bf Network theory}.
Graph theory has strong roots in computer science, where networks or geometric meshes in computer graphics are 
considered which are in general very complex. Social networks \cite{Jackson}, biological networks \cite{BLM}, the web, 
or triangularizations of surfaces used to display computer generated images are examples. Many mathematicians look at graphs as 
analogues of Riemann surfaces, as one-dimensional objects therefore. Similarly like the GAGA principle seeing correspondences
between algebraic and analytic geometry, there is a principle which parallels the geometry of graphs and
the geometry of manifolds or varieties. It is also a place, where some Riemann themes come together:
{\bf Riemannian manifolds}, {\bf Riemann-Roch}, {\bf Riemann-Hurwitz}; studying the {\bf Riemann Zeta function}
on cyclic graphs leads naturally to {\bf Riemann sums}. \\

{\bf Dynamics}.
We have seen dynamical systems on graphs given by discrete partial differential equations. 
It is the dynamics when studying the evolution of functions or discrete differential forms on graphs.
This is heat or wave dynamics or quantum dynamics. 
Dynamics also appears is as {\bf geometric evolution equations}, which describe families of 
geometries on a graph. While the graph itself can be deformed by discrete time steps which preserve
continuity, one can also deform the exterior derivative. If this is done in a
way so that the spectrum of the Laplacian stays the same, we study symmetries
of the quantum mechanical system. 
Heat, wave, Laplace, Maxwell or Poisson equations are based
on the Laplacian $L$. Classically, Gauss defined the gravitational field $F$ as
the solution to the Poisson equation $LF = \sigma$, where $\sigma$ is the
mass density. The Newton equations describe particles $u(t)$
satisfying the equation $u''= L^{-1} \sigma$. Better is a
Vlasov description where particles and mass are on the same footing
and  where we look at deformations $q(t)$ which satisfy the differential equation 
$q''(x)=-L^{-1} q$, a differential equation which makes sense also in the 
discrete. \\

{\bf Algebra}.
There are other places where graphs can be seen from the dynamical system point of view.
It is always possible to see a graph as an orbit of a monoid action on a finite set.
Sometimes this is elegant. For example, look at the graph on $Z_p$
generated by the maps $T(x)=3x+1$ and $S(x)=2x$ or the system generated
by a quadratic map $T(x)=x^2+c$.
This orbital network construction produces deterministic realistic networks.
These graphs relate to {\bf automata}, edge colored directed graphs with possible self loops
and multiple loops and encode a  monoid acting on the finite set.
\cite{Steinberg} who has developed a theory of finite transformation monoids calls
a graph generated by $T$ an {\bf orbital digraph} which prompted us to call the graph
{\bf orbital networks}.
That the subject has some number theoretical flavour has been indicated in \cite{KnillOrbital1,KnillOrbital2,KnillOrbital3}.
It shows that elementary number theory matters when trying to understand connectivity properties of
the graphs. We would not be surprised to see many other connections. \\

{\bf Cup product}. Exactly in the same way as in the continuum, one can define on 
the space $\Omega$ of antisymmetric functions a product $f \cup g$ which maps $\Omega_p \times \Omega_q$ to 
$\Omega_{p+q}$ and which has the property that the product of two cocycles is a cocycle extending so to a 
product $H^k \times H^l \to H^{k+l}$. The definition is the same: assume $f \in \Omega^p$ and $g \in \Omega^q$, 
then define $f \cup g(x_0,...,x_{p+q}) = f(x_0,...x_p) g(x_p,...x_{p+q})$. Now check that 
$d(f \cup g) = df \cup g + (-1)^p f \cup dg$. 
For example, if $p=q=1$, then $h(x_0,x_1,x_2) = (f \cup g)(x_0,x_1,x_2) = f(x_0,x_1) g(x_1,x_2)$ and 
$dh(x_0,x_1,x_2) = f(x_1,x_2) g(x_2,x_3) - f(x_0,x_2) g(x_2,x_3) + f(x_0,x_1) g(x_1,x_3) - f(x_0,x_1) g(x_1,x_2)$
which agrees with $df \cup g - f \cup dg 
= (f(x_1,x_2) - f(x_0,x_2) + f(x_0,x_1))  g(x_2,x_3) - f(x_0,x_1) ( g(x_1,x_2) - g(x_1,x_3) + g(x_2,x_3) )$. 
We see that if $df=0$ and $dg=0$, then $d(f \cup g)=0$. 
Lets take the example, where $G=(V,E)$ is the octahedron, where $V=\{a,b,c,d,p,n\}$ and $E$ consists of 12 edges. 
Define first a $1$-form $f$ which is equal to $1$ 
on the equator $a \to b \to c \to d \to a$ and zero everywhere else. Assume that the remaining edges are
oriented so that all point to the north or south pole $n$ rsp $p$. The second $1$-form $g$ is now chosen to 
be $1$ on each of the remaining $8$ edges and zero on the equator. We can fix an orientation now on the triangles 
$(x,y,z)$ so that $f \cup g(x,y,z) =1$ for every triangle. This is a nonvanishing $2$-form is a discrete area form
and showing that ${\rm cup}(G)=2$. Of course, one has also ${\rm tcap}(G)=2$ since we can find two
in $G$ contractible sets which cover $G$. The function $h$ given by $h(p)=0,h(a)=1,h(b)=2,h(c)=3,h(d)=4,h(n)=5$
has exactly two critical points $\{n,p\}$ because $p$ is the minimum with $i_h(x) = \chi(S^-(p)) = 1-\chi(\emptyset)=1$
and $n$ is the maximum with $i_h(n) = 1-\chi(S^-(n)) = 1-\chi(\{a,b,c,d \}) = 1-\chi(C_4) = 1$. This is an 
example, where both Ljusternik-Schnirelman inequalities are equalities: 
${\rm cup}(G) =  {\rm tcap}(G) = {\rm crit}(G)$.  \\

{\bf Calculus of variations}. One can study various variational problems for graphs. 
The most important example is certainly the Euler characteristic $\chi(G)$. It actually 
can be seen as a discrete Hilbert action \cite{eveneuler}. Network scientists look at
the characteristic length $\mu(G)$ or the mean clustering $\nu(G)$.
One can also look at the dimension $\iota(G)$, edge density $\epsilon(G)$, scale measure $\sigma(G)$,
or  spectral complexity $\xi(G)$ of a graph. 
There is a lot of geometry involved. The cluster coefficient $C(x)$ for example is closely related
to characteristic length. If $\mu$ denotes the mean distance and $\nu$ is the mean clustering coefficient and
$\eta$ is the average of scalar curvature $S(x)$, a formula $\mu \sim 1+\log(\epsilon)/\log(\eta)$
of Watts and Strogatz \cite{NewmanStrogatzWatts} relates $\mu$ with the edge density 
$\epsilon$ and average scalar curvature
$\eta$ telling that large curvature implies small average length, a fact familiar for spheres.
We often see statistical relations $\mu \sim \log(1/\nu)$ holds for random or
deterministic constructed networks, indicating that small clustering is often associated to large
characteristic lengths. Clustering $\nu$, edge density $\epsilon$
and curvature average $\eta$ therefore can relate with average length $\mu$ on a statistical level. \\

{\bf Pedagogy}.
Graphs are intuitive objects. We know them from networks like street maps
or subway networks. While calculus needs some training, especially to in higher dimensions,
the fundamental theorem of calculus in any dimension can be formulated and proven in a couple of minutes within
graph theory. Without any technicalities, it illustrates the main idea of Stokes theorem. While differential
forms need some time to be mastered by students learning analysis, functions on simplices of graphs are
very concrete in the sense that we always can get concrete results. A couple of lines in a computer algra system
provides code, which universally works in principle for any finite simple graph. Procedures spitting out a
basis for the cohomology groups for example is convenient since we just have to compute the eigenvalues and
eigenvectors of a concrete matrix $L$ and take the eigenvectors belonging to the eigenvalue $0$. 
The subject therefore has {\bf pedagogical merit}, also in physics, where wave and heat evolution on geometric
spaces can be visualized. If the spaces are discretizations of manifolds, the discretization provides a 
numerical scheme. Of course, the classical numerical methods to solve PDEs can do this often more 
effectively but then, one has in each case to deveop code for specific situations. Evolving waves in two 
and three dimensions effectively for example needs completely different techniques. The graph case is 
simple because the a dozen lines can deal with an arbitrary graph, as long as the machine can handle the 
matrices involved. For the computation of cohomology for example, we can first homotopically simplify the
graph as good as possible, then find a good \v{C}ech cover and then deal with a much smaller nerve graph. \\

{\bf Quantum calculus}.
Graph theory leads to relations between discrete mathematics, analysis and discrete differential geometry, 
analysis, topology and algebra. In a calculus setting, it leads to calculus without limits which is 
also called {\bf quantum calculus} because of commutation relations $[Q_e,P_e] = i$ 
which position $Q_e = q(x) \sigma^*$ and momentum operators $P_e f(x) = i [\sigma,f]= [f(x+e)-f(x)] \sigma$
attached to an edge $e=(x,y)$, where $q(y)=1$ and $q(x)=0$. For every path in the graph 
$\gamma$ we can define $Q_{\gamma}$ and have a {\bf linear approximation formula} $f(z) = f(x) + D_e f(x)$ for neighboring points $x,z$
and a {\bf Tylor formula} $x,z$ if we chose a path from $x$ to $z$ and have 
$f(z) = f(x) + D_e f(x) + D_{e_1} D_{e_2} f(x)/2! + ... = \exp(-i D) f$. This is analogue to the 
continuum where the Taylor equation means $\exp(-i t P) f(x) = f(x+t)$ if $P = i \partial_x$ is the momentum operator. 
Discrete Taylor formulas have already been known to Newton and his contemporary Gregory and are now mostly 
known in the context of numerical analysis. 
One can see the analysis of graphs as a ``quantization" since position and momentum operators stop to commute if 
space is discretized.  \\

{\bf Illustrations}.
 Graph theory allows to illustrate how different fields of mathematics
like algebra, calculus, analysis, topology, differential geometry, number theory or algebraic geometry can interplay.
Graph theory is also an applied topic in mathematics with strong roots in computer science, where networks or
geometric meshes in computer graphics are considered which are in general very complex. Social networks,
biological networks, the web, or triangularizations of surfaces used to display computer generated images
are examples.  Many mathematicians look at graphs as analogues of Riemann surfaces, as one
dimensional objects therefore. Similarly like the GAGA principle seeing correspondences
between algebraic and analytic geometry, there is principle which sees graphs and one
dimensional possibly complex curves and varieties as on the same footing.
But there is more to it: we can treat in many respects graphs like Riemannian manifolds or varieties.
It is also a place, where Riemann themes come together:
{\bf Riemannian manifolds}, {\bf Riemann-Roch}, {\bf Riemann-Hurwitz}; studying the {\bf Riemann Zeta function}
on cyclic graphs leads naturally to {\bf Riemann sums}. Which is one reason why I like the topic because
Riemann is one of my favorate mathematicians. \\

{\bf The point of view}.
Most mathematicians in geometric graph theory consider graphs as objects of
{\bf one-dimensional nature}: edges are one-dimensional arcs which 
connect the zero-dimensional vertices. Indeed, many results from algebraic curves parallel in graph 
theory. An example is the Riemann-Roch theorem \cite{BakerNorine2007} or the
Nowakowski-Rival fixed point theorem of 1979 which is a special case of a discrete Brouwer
fixed point theorem \cite{brouwergraph}. The Riemann-Roch theorem is of more {\bf algebro-geometric} nature, 
fixed point theorems are by nature interesting in {\bf dynamical systems theory} or {\bf game theory}.
In the last couple of years, the point of view started to shift and geometric
questions which are traditionally asked for higher dimensional geometric objects
started to pop up in graph theory. An example is {\bf discrete Morse theory} and
{\bf discrete differential geometry}. Morse theory deals studies critical points of
functions on a graph and relates them to geometry. Morse inequalities relate cohomology with 
the number of critical points of certain type, relating so algebra with analysis. 
An other example is {\bf category theory} in topology, which LS-category is a 
fascinating and central topic in topology because it links analysis, topology and algebra. 
An algebraic notion of how "rich" a space is the {\bf cup length}. It is defined in the cohomology ring of a 
graph. The topological notion is called LS-category which roughly tells with how many 
contractible sets one can patch up a space. The third notion is analytic and gives
the minimal number of critical points which a function can have on a graph. \\

{\bf Higher dimensional structures}.
The fact that notions of curvature, homotopy, degree or index \cite{brouwergraph}, or 
critical points can be defined for graphs so that classical results hold also in the discrete is a 
strong indication that there is more higher dimensional structure on a graph than anticipated. 
Until recently, this was only studied for smaller dimensions. For triangularizations of 
two-dimensional surfaces for example, there is a simple notion of curvature which
involves the degree of the graph. This has been known 
for a while \cite{princetonguide,Higuchi,Gromov87} but is probably much older in special
cases. The triangular lattice with hexagonal symmetry where every node has uniform 
degree $6$ for example is the prototype of a ``flat" geometry. 
The Gromov type curvature for planar graphs might have proposed in \cite{Stone}.
This curvature was generalized to arbitrary finite simple graphs \cite{cherngaussbonnet}. As indicated 
in \cite{elemente11}, first order notions of curvature can be refined. One can slso
average first order notions over smaller neighborhoods to get more refined notions. New definitions of
curvature based on probabilistic Markov chain concepts were introduced in \cite{Ollivier}
and pursued further in work like \cite{LinLuYau}. The Olivier curvature has proven to 
be fruitful and is under heavy investigation. 
It looks promising to get results close to classical differential 
geometry like that. The proposal under consideration uses other curvatures and plans to 
investigate the relations between different curvatures.
We have seen that virtually any 
topic in differential geometry can be studied in the discrete. Gauss-Bonnet,
Riemann-Hurwitz, fixed point theorems, Poincar\'e-Hopf, differential equations, isospectral graphs for 
the Dirac operator, integrable systems obtained by doing isospectral deformations, 
zeta functions for graphs or a McKean-Singer result using the heat kernel in graph theory. \\

{\bf Networking with structures}.
Since finite simple graphs are networks, there are connections with network geometry, 
which is itself close to fields like combinatorics, statistical mechanics or complexity theories.
What is particularly exciting is that the topic allows with modest techniques to illustrate how different 
fields of mathematics can overlap and play together. Gauss-Bonnet illustrates how the
metric property of curvature which is not topological is related to the topological property Euler characteristic.
The probabilistic relation with Poincar\'e-Hopf shows a integral geometric statistical angle. 
Integral calculus has penetrated classical differential geometry a long time ago and there are generations
of mathematicions like Blaschke-Chern-Banchoff which illustrate this. 
Ljusternik-Schnirelmann shows that algebraic, topological and analytic notions work together. 
The cohomologically defined cup product, the topologically defined cup length using covers and 
the analytically defined minimal number of critical points. 
In geometry, analysis deals with spaces of functions on a geometric space and is a theory called functional analysis.
Calculus and differetial topology deals with especially critical points. 
In the discrete as well as in more general situations like metric spaces, 
a notion homotopy leads immediately leads to a notion of critical points as points where
the part of small spheres for which the function value is smaller are not contractible. \\

{\bf The Dirac operator}.
An other important connection is through the Dirac operator $D$ which is more fundamental
than the Laplacian $L=D^2$ because its symmetry group is much larger than the one of the Laplacian. 
It allows to be deformed in an isospectral way without that the Laplacian is affected. This means
that these internal symmetries do not affect classical quantum mechanics but do affect the wave
equation in a particular way. Space expands. We can see that in the D'Alembert solution of the wave
equation where $D$ appears in the solution. In the expanded space, we have to prepare larger initial
velocities to get to the same solution. One can see the expansion also in the Connes reformulation of 
Riemannian geometry. 
Also in the Riemannian manifold case we have isospectral deformations of the Dirac operator $d+d^*$
and so a deformation of the exterior derivative which allows geometric evolutions to happen.
Cohomology of course does not change.  This is natural since these deformations are fundamental 
symmetries of the underlying geometry. Both in the continuum as well as in the discrete it produces 
an expansion of the metric with inflationary start. \\

{\bf What's next?} As mentioned in \cite{knillmckeansinger}, one goal is to have
a translation of Atiah-Singer which can be taught without much technical overhead. Already 
Gauss-Bonnet can be seen as a very special case of an index theorem, where the analytic index of the Dirac matrix
$D: \Omega_{even} \to \Omega_{odd}$ is the cohomological Euler characteristic and the
average curvature can be seen as a topological index. While Gauss-Bonnet-Chern is an almost
trivial case in the discrete, the next step is Hirzebruch-Riemann-Roch. It has been stated
as an open problem in \cite{BakerNorine2007} and it is not expected to be so simple also in the
discrete because more structure is needed like higher dimensional versions of divisors and an
adaptation of sheaf cohomology. Since the discrete case of classical theorems like Gauss-Bonnet
or Poincar\'e-Hopf is so short that full proofs could be presented in the first part of this 20
minutes Baltimore talk and where concrete and short computer algebra 
implementations \cite{KnillWolframDemo2,KnillWolframDemo1} containing detailed code doing this
this for arbitrary graphs, I would not be surprised
to finally see and understand Riemann-Roch in a discrete setup. This means that we can program it for
an arbitrary graph and divisor and get two quantities from it, which the theorem shows to be equal. 
While I had to learn as a student the Patodi proof of Gauss-Bonnet-Chern inf\cite{Cycon}, 
I don't understand Riemann-Roch yet in higher dimensions. It is an other level of difficulty and there 
are so many different aspects of the theorem.  Thinking about the discrete could
be a way to learn it. Having mentioned positive curvature graphs, there is a whole more to be explored, if the 
positive curvature assumption is relaxed to average positive curvature assumptions for subgraphs. 
One combinatorial question in the positive curvature case is to make a list of all positive curvature 
graphs in dimensions $d$ and to get sphere theorems.  The most exciting problem by far is to 
understand the isospectral deformation of the Dirac operator. It is important because the deformation also
works for Riemannian manifolds.  The deformed exterior derivative defines a deformed geometry on the 
graph or manifold. The question is: what geometry are obtained asymptotically if the expansion is rescaled?

\bibliographystyle{plain}

\end{document}